%% file: AbaqusSFEM.tex
\documentclass[preprint,12pt]{elsarticle}

\usepackage{graphicx}
\usepackage{amssymb}
\usepackage{epstopdf}
\usepackage{listings}
\usepackage{lineno}

\biboptions{sort&compress}
\usepackage{etex} 
\usepackage{graphicx,subfigure,color,colordvi}
\usepackage{latexsym,amsbsy,epsfig,fancybox}
\usepackage{amstext}
\usepackage{amsfonts,textgreek}
\usepackage{mathrsfs}
\usepackage{mathtools}
\usepackage{stmaryrd,dsfont}
\usepackage{bbm}
\usepackage{url}
\usepackage{wrapfig,bm}
\usepackage{tikz}
\usepackage{pstricks}
\usepackage{caption}
\usepackage{psfrag}
\usepackage{upgreek}
\usepackage{isomath}
\usepackage{latexsym}
\usepackage{array}
\usepackage{multirow}
\usepackage{booktabs}
\usepackage{colortbl}
\usepackage{verbatim}
\usepackage{amsmath,amssymb,xspace}
\usepackage{pgfplots}
\usetikzlibrary{plotmarks}

\usepackage[margin=3cm]{geometry}

\newcommand{\eref}[1]{Equation~(\ref{#1})}

\newcommand{\fref}[1]{Figure~\ref{#1}}
\newcommand{\frefs}[1]{Figures~\ref{#1}}

\newcommand{\vm}[1]{\bm{\mathrm{#1}}} 
\newcommand{\bsym}[1]{\bm{#1}}

\renewcommand{\Re}{{\rm{I\!R}}}

\newcommand{\mat}[1]{\bm{\mathrm{#1}}} 
\newcommand{\transpose}{\mathrm{T}}

\newcommand{\smat}[2][ccccccccccccccccccccccccccccccccccccccccccccccccccc]{\left
[\begin{array}{#1}#2 \\ \end{array} \right]}

\tikzstyle{nicebox}=[draw=black!100, fill=white!10, rectangle, inner sep=4pt, inner ysep=16pt]
\tikzstyle{niceboxtitle}=[draw=black!100, fill=white, text=black, rectangle]

\usetikzlibrary{arrows,decorations.markings}

\definecolor{forestgreen}{RGB}{34, 139, 34}
\definecolor{lightgray}{gray}{0.92}

\newcommand{\xx}{\vm{x}}
\newcommand{\rmd}{\mathrm{d}}
\newcommand{\bigb}{\mat{B}}

\newcommand{\dd}{\mat{D}}
\newcommand{\cn}{\vm{n}}

\newcommand{\bq}{\mat{q}}
\newcommand{\uu}{\mat{u}}

\newcommand{\bveps}{\bsym{\varepsilon}}
\newcommand{\bvsig}{\bsym{\sigma}}

\journal{Elsevier}

\begin{document}

\begin{frontmatter}


\title{An Abaqus UEL implementation of the smoothed finite element method}




\author{Pramod Y. Kumbhar\fnref{label1}}
\author{A. Francis\fnref{label1}}
\author{N. Swaminathan}
\author{R. K. Annabattula\corref{cor1}}
\ead{ratna@iitm.ac.in}
\author{S. Natarajan\corref{cor1}}
\ead{snatarajan@iitm.ac.in}

\cortext[cor1]{Corresponding author}
\fntext[label1]{Equal contribution}
\address{Integrated Modelling and Simulation Lab, Department of Mechanical Engineering, Indian Institute of Technology Madras, Chennai 600036, India.}

\begin{abstract}

In this paper, we discuss the implementation of a cell based smoothed finite element method (CSFEM) within the commercial finite element software Abaqus. The salient feature of the CSFEM is that it does not require an explicit form of the derivative of the shape functions and there is no isoparametric mapping. This implementation is accomplished by employing the user element subroutine (UEL) feature of the software. The details on the input data format together with the proposed user element subroutine, which forms the core of the finite element analysis are given. A few benchmark problems from linear elastostatics in both two and three dimensions are solved to validate the proposed implementation. The developed UELs and the associated input files can be downloaded from 
\href{Github repository link}{https://github.com/nsundar/SFEM\_in\_Abaqus}.
\end{abstract}


\end{frontmatter}


\section{Introduction}
\label{intro}
The finite element method (FEM) is probably the most popular technique to numerically solve partial and ordinary differential equations in science and engineering. This popularity is due to the flexibility the method provides in treating complex geometries and boundary conditions. The method involves discretizing the domain into several non-overlapping regions called the elements that are connected together by a topological map called a mesh. Conventional finite elements are generally restricted to simplex elements, viz., triangles and quadrilaterals in two dimensions and tetrahedra and polyhedra in three dimensions with a local polynomial representation used for the unknown field variables within each element. Despite its popularity, the FEM with simplex elements faces the following challenges: (a) the accuracy of the solution depends on the quality of the element\cite{leebathe1993}, i.e., the results are sensitive to mesh distortion and (b) requires sophisticated discretization techniques to generate high quality meshes and to capture topological changes.

To alleviate some of the aforementioned challenges, recent research has been focused on alternative approaches, viz., meshfree methods~\cite{Belytschkolu1994,zhuzhang2007,aluru1999,khosravifardhematiyan2017,Nguyen2017}, partition of unity methods (PUM)~\cite{babuskabanerjee2004,belytschofries2010,bordasmoran2006,tanakasuzuki2016,Tanaka2017}, strain smoothing techniques~\cite{liudai2007,shengmin2003} and polygonal finite element method (PFEM)~\cite{francisa.ortiz-bernardin2017,sukumartabarraei2004,Leekim2016,veigabrezzi2013}. The meshfree methods does not have a priori topology which makes it suitable for large deformation problems and moving boundary problems~\cite{nguyenrabczuk2008,rabczukzi2010,rabczukblytschko2004}. The PUMs relaxes the meshing constraint by allowing the geometric features to be represented independent of the underlying FE discretization, whilst PFEM relaxes the meshing constraint by allowing elements to take arbitrary shapes~\cite{talischipaulino2012}. 

On the other front, the strain smoothing technique improves the quality of the FE solution by computing a modified strain field which is a weighted average of the compatible FE strain field~\cite{liudai2007}. Inspired by the application of the strain smoothing to meshfree methods~\cite{chenwu2001}, Liu and co-workers~\cite{liudai2007} introduced the technique within the FEM and coined the method, `smoothed finite element method (SFEM)'. In this method, the modified strain is written over the smoothing domain. The method can also be seen as dividing the domain into smoothing cells, however, from a practical point of view,  it is computationally less expensive when using the existing mesh data structure to generate the smoothing domain. It has been shown theoretically and numerically that the SFEM is more accurate than the conventional FEM~\cite{liunguyen2010}.

Liu \textit{et al.}~\cite{liunguyen2007} formulated a series of SFEM models based on the choice of smoothing domain, viz., cell-based SFEM (CSFEM)~\cite{liudai2007}, node based SFEM (NSFEM)~\cite{liuthoi2009}, edge-based SFEM (ESFEM)~\cite{liuthoi2009a}, face-based SFEM (FSFEM)~\cite{thoiliu2009}, $\alpha$FEM~\cite{liunguyen2008} and $\beta$FEM~\cite{zengliu2016}. 
\begin{figure}[ht]
\centering
\includegraphics[scale=0.6]{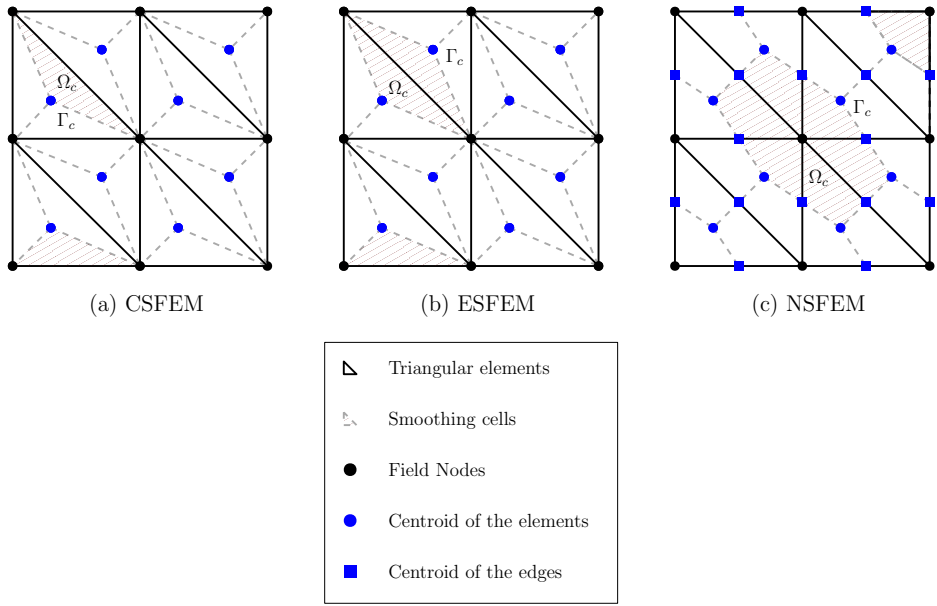}
\caption{Schematic representation of the smoothing cells $\Omega_c$ over triangular elements for the cell-based, edge-based and node-based smoothing technique. The `shaded' region is the smoothing cell, over which the smoothed strains are computed from the compatible strain field.}
\label{fig:sfemillustration}
\end{figure}
\fref{fig:sfemillustration} shows a schematic various strain smoothing techniques in two dimensions. The common feature of these SFEM models are that all these models use finite element meshes with linear interpolants, except for the CSFEM. The choice of smoothing technique controls the assumed strain field and thus yielding accurate and stable results~\cite{bordasnatarajan2010}. Recently, the CSFEM has been applied to arbitrary polygons and polyhedra~\cite{francisa.ortiz-bernardin2017}. An interesting feature of using polygons and polyhedra together with the strain smoothing technique is the following. The nodal strain is computed through the divergence of a spatial average of the standard local strain field. This aspect of the technique allows reduces the integration space by an order. Consequently, the stiffness matrix is computed on the boundaries of the smoothing elements instead of the whole  domain.

Recently, the SFEM has also been incorporated into the PUM~\cite{bordasrabczuk2010,bordasnatarajan2011,chenrabczuk2012,wuliu2016} and PFEM~\cite{dailiu2007}. Liu and co-workers~\cite{liubui2012,nguyenliu2013} have applied the ESFEM to solve problems involving strong discontinuities with varying orders of singularity. While these methods exist as stand-alone codes, they have not yet formed a part of a commercial software and this aspect restricts the use of the method to a specific community or laboratory. Therefore, in this work we implement the technique within the commercial finite element software Abaqus and provide the details so that other researchers can easily extend this work for more specific or general cases. The procedures that can be used by practicing engineers who are familiar with Abaqus can get benefited with the capabilities of the SFEM. For the sake of simplicity, we choose the most commonly used smoothing technique, the cell based SFEM (CSFEM). 

Although the method is general the implementation aspects of the SFEM are applied to two and three dimensional elasticity to demonstrate the technique. Particular emphasis is placed on the input data format and the subroutines that implement the smoothing technique. 
The method will be implemented over polytopes 
to demonstrate the accuracy and the robustness of the smoothing technique when compared to the conventional FEM. 

The remainder of the manuscript is organized as follows. The next section presents the governing equations for elastostatics and the corresponding weak form. Section \ref{sfembasics} presents an overview of the cell based smoothing technique over arbitrary polygons and polyhedra. The Abaqus UEL implementation is discussed in Section \ref{abaqimplementation}. Section \ref{numersection} presents a detailed convergence study by comparing the accuracy and the convergence rates with the conventional FEM for a few benchmark problems in two and three dimensional linear elastostatics, followed by concluding remarks in the last section.

\section{Governing equations and weak form for linear elastostatics}
\label{governeqn}
In this section, the governing equations for linear elastostatics are first presented, followed by the corresponding weak form. Let us consider an isotropic linear elastic body in $d$ ($\{=2,3\}$) dimensional Euclidean space, $\Re^d$, whose material point is given by $\xx = \sum \xx_I \mathbf{e}_I$, where $\mathbf{e}_I$ are the vectors of a chosen basis. Let $\Omega \subset \Re^d$ represent the body with domain boundary $\Gamma \equiv \partial \Omega$ consisting of Dirichlet boundary $\Gamma_u$ and Neumann boundary $\Gamma_t$, such that $\Gamma=\Gamma_u \cup \Gamma_t$, $\emptyset=\Gamma_u\cap\Gamma_t$ and the outward normal to $\Gamma$ is $\vm{n}$. Then, the strong form of the boundary value problem is: given the body force $\vm{b}$, the Cauchy stress tensor $\bvsig$ and the prescribed displacement $\bar{\vm{u}}$ and the traction $\bar{\vm{t}}$, find the displacement field $\vm{u} : \bar{\Omega} \rightarrow \Re^d$ such that
%
%
%
\begin{align}
\quad \bsym{\nabla} \cdot \bsym{\sigma} + \vm{b} &= \textbf{0} ~ \quad \,\textup{in} \quad \, \Omega, \nonumber \\
\quad \vm{u} &= \bar{\vm{u}} ~ \quad \, \textup{on} \quad \, \Gamma_u, \nonumber \\
\quad \bsym{\sigma} \cdot \vm{n}
&= \bar{\vm{t}} ~ \quad \,\textup{on} \quad \, \Gamma_t,
\label{eqn:problem_strong_form}
\end{align}

Let $\mathscr{U}(\Omega) = \left\{ \vm{u}: \Omega \rightarrow \Re^d | u_I \in H^{1}(\Omega), I=1,\cdots,d , \ \vm{u} = \bar{\vm{u}} \ \textrm{on } \Gamma_u  \right\}$ be the displacement trial and $\mathscr{V}(\Omega) = \left\{ \vm{v}: \Omega \rightarrow \Re^d | v_I \in H^{1}(\Omega), I=1,\cdots,d, \ \vm{v} = \vm{0} \ \textrm{on } \Gamma_u \right\}$ be the test function spaces, where $H^{1}$ denotes the Hilbert-Sobolev first order space. The resulting weak form can be written as: find $\vm{u} \in \mathscr{U}$  such that:
\begin{equation}
\quad a(\vm{u},\vm{v}) = \ell(\vm{v}) ~\forall \vm{v} \in \mathscr{V}
\label{eqn:weakform}
\end{equation}
where the bilinear and linear forms are defined as follows:
\begin{align}
\quad a(\vm{u},\vm{v})&=\int_{\Omega}\bsym{\sigma}(\vm{u}):\bsym{\varepsilon}(\vm{v})\,\rmd \Omega, \nonumber \\
\quad \ell(\vm{v})&= \int_{\Omega}\vm{b}\cdot\vm{v}\,\rmd \Omega + \int_{\Gamma_t}\bar{\vm{t}}\cdot\vm{v}\,\rmd \Gamma,
\label{eqn:bilinearlinearform}
\end{align}
where $\bsym{\varepsilon}=\frac{1}{2} \left[ \nabla \uu + \nabla \uu^{\rm T} \right]$ is the small strain tensor. In the following, it is assumed that the arbitrary domain is partitioned into $nel$ finite elements defined as $\bar{w}_I$ such that $\Omega \equiv {\rm{int}(\cup_{I=1}^{nel}} \bar{w}_I)$ and $\bar{w}_I \cap \bar{w}_J = \emptyset, \forall I \neq J$. For the discretization of the weak form (see \eref{eqn:weakform}), let the set $\mathscr{U}^h \subset H^1(\Omega)$ consist of polynomial interpolation functions of the following form:
\begin{align}
\vm{u}^h=\sum_{e=1}^{nel} \phi_e \vm{u}_e \nonumber \\
\vm{v}^h=\sum_{e=1}^{nel} \phi_e \vm{v}_e
\label{eqn:dispfeapprox}
\end{align}
where $\phi_e$ denote the finite element shape functions. Upon substituting \eref{eqn:dispfeapprox} into~\eref{eqn:weakform}, the following discrete weak form is obtained, which consists of finding $\vm{u}^h \in \mathscr{U}^h$ such that
\begin{equation}\label{eq:discweakform}
 \quad a(\vm{u}^h,\vm{v}^h)= \ell(\vm{v}^h) ~ \forall \vm{v}^h \in \mathscr{V}^h,
\end{equation}
which leads to the following system of linear equations:

\begin{subequations}
\begin{align}
\begin{split}
\mat{K}\mat{u}&=\mat{f}, \label{eq:weakform_disc} 
\end{split} \\
\begin{split}
\hspace{-0.6cm}\mat{K}&=\sum_h\mat{K}^h=\sum_h\int_{\Omega^h}\mat{B}^\transpose\dd \mat{B}\,\rmd \Omega, \\
\end{split} \\
\begin{split}
\mat{f} &= \sum_h \mat{f}^h=\sum_h\left(\int_{\Omega^h}{\boldsymbol{\phi}}^\transpose\vm{b}\,\rmd \Omega + \int_{\Gamma_t^h}{\boldsymbol{\phi}}^\transpose\bar{\vm{t}}\,\rmd \Gamma\right), 
\end{split}
\end{align}
\end{subequations}
where $\mat{K}$ is the global stiffness matrix, $\mat{f}$ is the global nodal force vector, $\mat{u}$ the global vector of nodal displacements, ${\boldsymbol{\phi}}$ is the matrix of finite element shape functions, $\dd$ is the constitutive relation matrix for an isotropic linear elastic material, and $\mat{B}=\bsym{\nabla} {\boldsymbol{\phi}}$ is the strain-displacement matrix that is computed using the derivatives of the shape functions. In the conventional FEM, the shape functions, $\phi_e$ are defined over the standard element and by invoking the isoparametric formulation, the integrals in~\eref{eq:weakform_disc} are evaluated using numerical integration. 

In this paper, we discretize the domain with arbitrary polygons and polyhedra in two and three dimensions, respectively. While there are different ways to represent the shape functions over arbitrary polytopes~\cite{sukumarmalsch2006}. we choose the Wachspress interpolants to describe the unknown fields~\cite{Wachspress1975}. These functions are rational polynomials and the construction of the coordinates is as follows: Let $P \subset \Re^3$ be a simple convex polyhedron with facets $F$ and vertices $V$. For each facet $f \in F$, let $\cn_f$ be the unit outward normal and for any $\xx \in P$, let $h_f(\xx)$ denote the perpendicular distance of $\xx$ to $f$, which is given by
\begin{equation}
h_f(\xx) = (\mathbf{v}-\xx) \cdot \cn_f
\end{equation}
for any vertex $\mathbf{v} \in V$ that belongs to $f$. For each vertex $\mathbf{v} \in V$, let $f_1,f_2,f_3$ be the three faces incident to $\mathbf{v}$ and for $\xx \in P$, let
\begin{equation}
w_{\mathbf{v}}(\xx) = \frac{ \mathrm{det}(\mat{n}_{f_1},\mat{n}_{f_2},\mat{n}_{f_3})}{ h_{f_1}(\mat{x}) h_{f_2}(\mat{x}) h_{f_3}(\mat{x})}
\end{equation}
with a condition that the ordering of $f_1,f_2,f_3$ be anticlockwise around the vertex $\mathbf{v}$ when seen from outside $P$. Then the barycentric coordinates for $\xx \in P$ is given by~\cite{warren2003,warrenschaefer2007}:
\begin{equation}
\phi_{\mathbf{v}}(\xx) = \frac{ w_{\mathbf{v}}(\xx)}{ \sum\limits_{\mathbf{u} \in V} w_{\mathbf{u}}(\xx)}.
\label{eqn:wachs3d}
\end{equation}
\begin{figure}[htbp]
\centering
\includegraphics[scale=0.8]{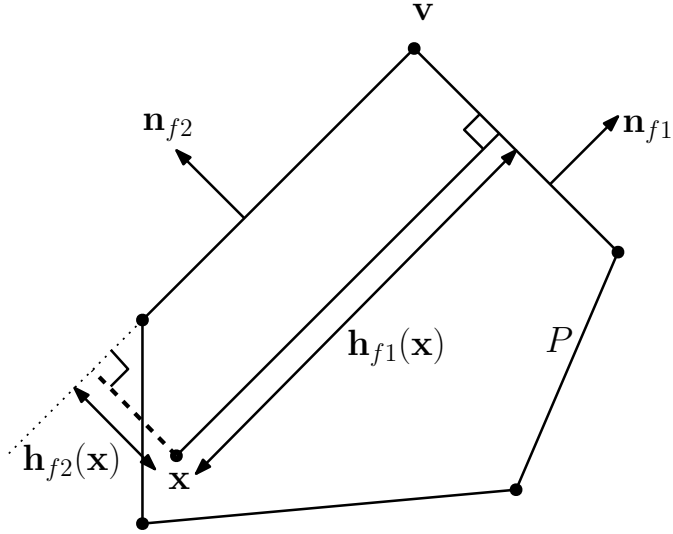}
\caption{Barycentric coordinates: Wachspress basis function}
\label{fig:bary}
\end{figure}
The Wachspress shape functions are the lowest order shape functions that satisfy boundedness, linearity and linear consistency on convex polytopes~\cite{warren2003,warrenschaefer2007}

\section{Overview of the smoothed finite element method}
\label{sfembasics}
In the strain smoothing technique, the compatible strains in the FEM are no longer used to compute the terms in the stiffness matrix. Instead, the integration related to the computation of stiffness matrix is based on the smoothing domain and hence no longer associated with the elements. 

Within the SFEM framework, the discrete modified strain field $\tilde{\varepsilon }_{ij}^h$ that yields the modified strain-displacement matrix $(\tilde{\mat{B}})$ which is then used to build the stiffness matrix is related to the compatible strain field ${\varepsilon }_{ij}^h$ by:
\begin{equation}
\tilde {\varepsilon }_{ij}^h (\xx)=\int_{\Omega_C^h}
{\varepsilon _{ij}^h (\xx) ~\bq(\xx) \rmd \Omega }
\label{eqn:epsilonvar}
\end{equation}
where $\bq(\xx)$ is the smoothing function and $\Omega_C^h$ is the sub-cell. On writing \eref{eqn:epsilonvar} at the basis functions derivative level and invoking Gauss-Ostrogradsky theorem, we get:
\begin{equation}
\int_{\Omega_C^h}\phi_{I,x}~\bq(\vm{x})\, \rmd \Omega  =  \int_{\Gamma_C^h}\phi_I~ \bq(\vm{x})n_j\, \rmd \Gamma - \int_{\Omega_C^h}\phi_I ~\bq_{,x}(\vm{x})\, \rmd \Omega
\label{eq:divconsistency}
\end{equation}
There are several choices for this $\bq(\xx)$. For simplicity, $\bq(\xx)$ satisfies the following properties:
\begin{equation}
 \bq \geq0    \hspace{0.2cm} \&    \int\limits_{\Omega^h}\bq(\xx) \rmd \Omega = 1
 \label{eqn:smoothcond}
\end{equation}
In this work, the weighting function, $\bq(\xx)$ is chosen to be a constant, i.e, $\bq(\xx) = 1$ for simplicity. For $\bq(\xx)=1$, the modified strains computed using \eref{eq:divconsistency} are constant within the smoothing domain. The Wachspress interpolants are used in this formulation.
%
%
%
%

\begin{figure}[h!]
\centering
\includegraphics[scale=0.75]{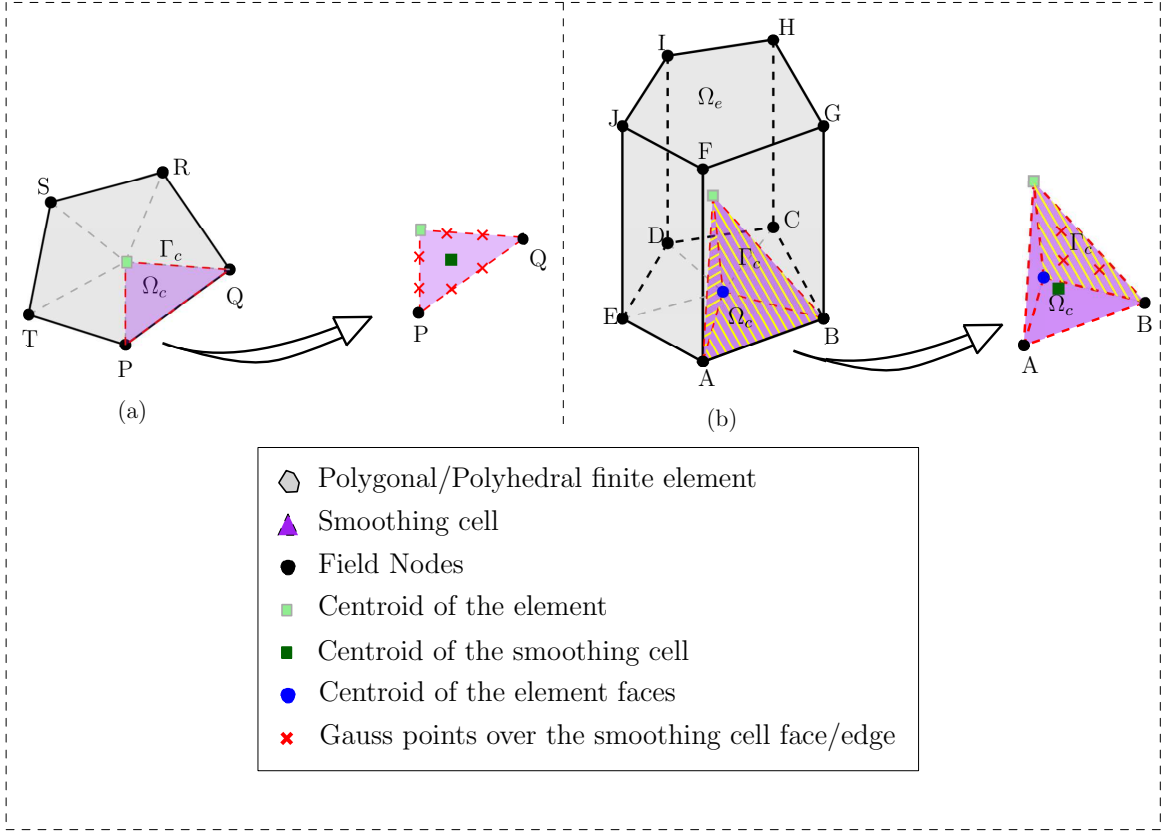}
\caption{Schematic representation of the smoothing cells over two and three dimensional polygonal element. In this study, we employ triangular and tetrahedral subcells, in two and in three dimensions, respectively.}

\label{fig:csfemsc}
\end{figure}

In the CSFEM, the local smoothing domain is constructed by sub-dividing the physical element into number of subcells. A subcell for an element is constructed by identifying the centroid\footnote{Here we use the geometric centre of the element.} of the element and connecting the vertices of the polygon/polyhedra with the geometric center of the element. \frefs{fig:csfemsc} 
shows the subcells within a polygon and a polyhedra.

Over each subcell, \eref{eq:divconsistency} is evaluated numerically to get the modified strain-displacement matrix. The stiffness matrix computed as for the CSFEM is:
\begin{equation}
\tilde{\mat{K}} = \sum\limits_{I=1}^{nel} \tilde{\mat{K}^h} = \sum\limits_{I=1}^{nel} \int\limits_{\Omega^e} \tilde{\mat{B}}^\transpose\dd \tilde{\mat{B}}\,\rmd \Omega,
\label{eqn:kmatrix}
\end{equation}
where $\tilde{\mat{B}}$ is the modified strain-displacement matrix. For the chosen smoothing function, the derivative of the function vanishes within the subcell and one integration point (i.e centroid of the smoothing cell $\Omega_c$ denoted by $a$) is sufficient to compute the modified strain-displacement matrix and the corresponding smoothed elemental stiffness matrix. 

Upon substituting the smoothing function, $\bq(\xx) = 1$ in \eref{eq:divconsistency}, the terms in the modified strain-displacement matrix can be computed by:
\begin{equation}
\phi_{a,i} = \dfrac{1}{A_c} \int\limits_{\Gamma_C^h}\phi_a(\vm{x}) n_i\,\rmd \Gamma
\label{eq:intconstraints}
\end{equation}
where $i=x,y$ for two dimensions and $i=x,y,z$ for three dimensions, $n_i$ are the unit normals to the edge and face in two and three dimensions, respectively and $A_c$ is the area and volume of the subcell in two and three  dimensions respectively. The above equation for a polygon of $n$ sides, gives the following expression for the modified strain-displacement matrix evaluated at the interior Gauss point:
\begin{equation}
\tilde \bigb=\smat{\tilde {\bigb}_1 &
\tilde {\bigb}_2 & \cdots & \tilde {\bigb}_n} \label{eq:Btilde3dbis},
\end{equation}
where the nodal matrix is
\begin{equation}
\tilde{\bigb}_a =\frac{1}{A_c} \int_{\Gamma_C^h}
\smat{n_x(\xx) & 0 \\
 0 & n_y(\xx) \\
 n_y(\xx) & n_x(\xx) } \phi_a(\xx) \rmd \Gamma
\label{eqn:Btilde2d}
\end{equation}
for two dimensions and
\begin{equation}
\tilde{\bigb}_a =\frac{1}{A_c} \int_{\Gamma_C^h}
\smat{n_x(\xx) & 0 & 0\\
 0 & n_y(\xx) & 0\\
 0 & 0 & n_z(\xx) \\
 n_y(\xx) & n_x(\xx) & 0 \\
 0 & n_z(\xx) & n_y(\xx) \\
 n_z(\xx) & 0 & n_x(\xx) } \phi_a(\xx) \rmd \Gamma
\label{eqn:Btilde3d}
\end{equation}
for three dimensions. From the above equations, it is clear that only the shape functions are involved in calculating the terms in the stiffness matrix. To evaluate the surface integrals in \eref{eqn:Btilde3d}, we employ the conforming interpolant quadrature technique.

\section{Abaqus UEL implementation:}
\label{abaqimplementation}
In this section, we describe the major steps involved in the implementation of UEL for the CSFEM in Abaqus. We focus on the subroutines and the input files which go with it. 
Through UEL, the stiffness and the residual contributions of an element to the system as stated in the Abaqus manual~\cite{abaqusmanual} are provided. The Abaqus solver can be invoked by the command:
\begin{center}
\textit{abaqus job=\textless inp file name \textgreater  user=\textless fortran file name\textgreater}. 
\end{center}
%
%
%
%
\begin{figure}[h!]
\centering
\includegraphics[scale=0.7]{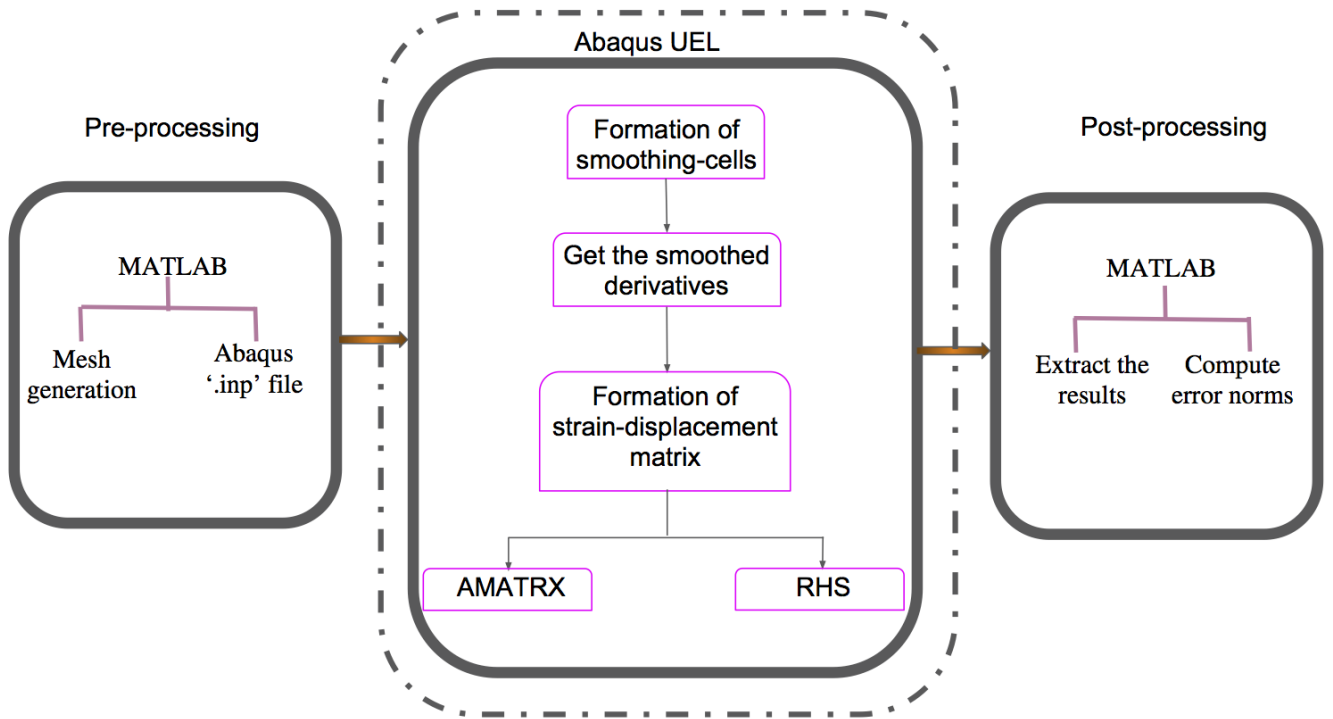}
\caption{Schematic illustrating the overall process involved in the implementation of the CSFEM within Abaqus.}
\label{fig:processchart}
\end{figure}
\fref{fig:processchart} depicts the three steps involved in the implementation of the smoothing technique within the Abaqus. The first stage is the pre-processing, where a MATLAB routine is used to generate the arbitrary polygonal and polyhedra meshes and the Abaqus input file. The second stage involves the core of the implementation, where the strain smoothing is implemented in UEL. As Abaqus CAE does not support post-processing of unconventional elements without significant modifications, a MATLAB script by Abaqus2matlab tool~\cite{papazafeiropouloscalvente2017} is used to extract results and compute necessary error norms.

\subsection{Structure of the input file}
The Abaqus input file contains the following information: nodal coordinates, element connectivity, material properties (Young's modulus and Poisson's ratio), boundary conditions and the solution procedure to be used. The boundary conditions and the solution procedures are specified following the standard procedure used in {\emph{Abaqus Standard}}. A new user element in Abaqus is defined through a key word {\emph{*User Element}}, the arguments to which are (in sequence):
\begin{itemize}
\item the number of nodes ($n$ for polygonal element having $n$ sides).
\item user assigned unique label starting with `U' followed by an integer. In the present implementation, the integer is tagged to the number of nodes of the element.
\item number of properties associated with the element (in this case $E$ and $\nu$).
\item the number of degrees of freedom per node (2 for 2D and 3 for 3D).
\end{itemize}

\begin{figure}
\centering
\includegraphics[scale=1]{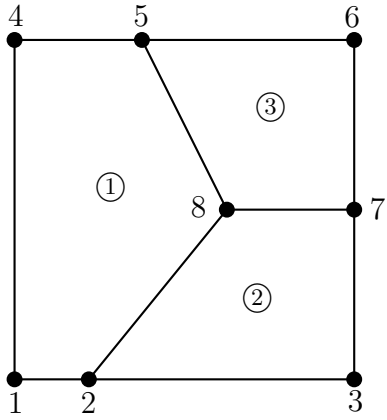}
\caption{A sample polygonal finite element mesh used for CSFEM. The CSFEM implementation requires only nodal coordinates and element connectivity.}
\label{suppl_file_mesh}
\end{figure}

In the case of CSFEM, arbitrary polygonal elements are used and hence the elements having same attributes (i.e., number of nodes) are grouped together with a unique label in the input file (see Listing~\ref{inp_file}). For a sample mesh shown in  \fref{suppl_file_mesh}, element 1 is tagged with label `U5' and elements 2 and 3 are grouped together with a label `U4'. The corresponding element connectivity is given after the keyword '{\emph{*Element}}'. 
\begin{center}
\begin{lstlisting}[caption=User element definition for arbitrary polygons (c.f. \fref{suppl_file_mesh}),label=inp_file]
*User element, nodes=5, type=U5, properties=2,coordinates=2
1,2
*Element, type=U5,ELSET=five
1,1,2,8,5,4
*UEL Property, ELSET=five
3.0e+07, 0.30
*User element, nodes=4, type=U4, properties=2,coordinates=2
1,2
*Element, type=U4,ELSET=four
2,2,3,7,8
3,8,7,6,5
*UEL Property, ELSET=four
3.0e+07, 0.30
\end{lstlisting}
\end{center}

\subsection{ User element details}
The UEL for implementing the CSFEM is written in FORTRAN90. The UEL starts with the name of the subroutine with the list of arguments in the parentheses~(see Listing \ref{subR}). The most important inputs to be provided through the UEL are RHS and AMATRX. Through RHS, the residual force contribution for each element is provided to the solver. AMATRX represents the element stiffness (or more generally the Jacobian) matrix as given by \eref{eqn:kmatrix}, where the strain-displacement matrix for CSFEM is given by~\eref{eq:Btilde3dbis}. The Abaqus solver constructs the global stiffness matrix by calling the UEL for each element  based on the connectivity information provided in the input file. 

     
\begin{center}
\begin{lstlisting}[caption=Snippet of the user element subroutine,label=subR]
     subroutine uel(rhs,amatrx,svars,energy,ndofel,nrhs,nsvars,
     1     props,nprops,coords,mcrd,nnode,u,du,v,a,jtype,time,dtime,
     2     kstep,kinc,jelem,params,ndload,jdltyp,adlmag,predef,
     3     npredf,lflags,mlvarx,ddlmag,mdload,pnewdt,jprops,njprop,
     4     period)
!     
      include 'aba_param.inc'
!

         dimension rhs(mlvarx,*),amatrx(ndofel,ndofel),
     1     svars(nsvars),energy(8),props(*),coords(mcrd,nnode),
     2     u(ndofel),du(mlvarx,*),v(ndofel),a(ndofel),time(2),
     3     params(3),jdltyp(mdload,*),adlmag(mdload,*),
     4     ddlmag(mdload,*),predef(2,npredf,nnode),lflags(*),
     5     jprops(*)

!     *variables for CSFEM*

!     * get the normals of the faces
      call getunitnorm(ncoord,...,g)

! loop over the faces to get the smoothed derivatives
! This is computed by employing Equation 17
! for two dimensions and Equation 18 for three dimensions

      call getsmoothderivative(ncoord,...,dNdz)

!     formation of stiffness(AMATRX) and residual(RHS) matrix

          end
                    
          subroutine getunitnorm(anodef,...,ag)
          ...
          end subroutine getunitnorm
          
          subroutine getsmoothderivative(ncoord,...,dNdz)
          ...
          end subroutine getsmoothderivative
\end{lstlisting}
\end{center}
\section{Numerical Examples}
\label{numersection}
The convergence properties and the accuracy of the CSFEM applied to arbitrary polytopes have been studied in detail in the literature~\cite{nguyenbordas2008,liudai2007,natarajanbordas2015}. The main objective of this section is to validate the implementation by solving a few benchmark problems in linear elastostatics. It is noted that readers can utilize the UEL subroutine to their specific applications with minimal modifications. For all the numerical examples, the domain is discretized with arbitrary polytopes. The two dimensional polygonal meshes are generated by using the built-in MATLAB function {\textit{voronoin}} and the MATLAB functions in Polytop~\cite{talischipaulino2012}. For three dimensional problems, the Polytop function is extended to generate the polyhedra meshes. 

%
%
%
%
%

For validation, the results from the CSFEM are compared with the corresponding regular FEM. The accuracy and the convergence rate are analyzed for different methods using the relative error norms in the displacement and the energy as given by\\
{\emph{Displacement norm}}:
\begin{equation}
|| \uu - \uu^h||_{\mathrm{L}^2(\Omega)} = \frac{ \sqrt{ \int\limits_{\Omega} (\uu-\uu^h)^{\rm T}  (\uu-\uu^h)~\rmd \Omega}}{\sqrt{ \int\limits_{\Omega} \uu ^{\rm T} \uu~\rmd \Omega}},
\end{equation}\\
{\emph{Energy norm}}:
\begin{equation}
|| \bveps - \bveps^h||_{\mathrm{H}^1(\Omega)} = \frac{ \sqrt{ \int\limits_{\Omega} (\bveps - \bveps^h)^{\rm T} \dd (\bveps - \bveps^h)~\rmd \Omega}}{\sqrt{ \int\limits_{\Omega} \bveps^{\rm T} \dd \bveps~\rmd \Omega}},
\end{equation}
where $\uu, \bveps$ are the analytical or reference solutions and $\uu^h, \bveps^h$ are the corresponding numerical solutions. The following convention is employed whilst discussing the results:
\begin{itemize}
\item PFEM-$n$d - polygonal finite element method. Here the numerical integration is performed by sub-triangulation method.
\item CSFEM-$n$d- cell based SFEM over arbitrary polytopes.
\end{itemize}
where $n=2,3$ represents the dimensionality of the problem.  


\subsection{Two dimensional Cantilever beam subjected to end shear}
\begin{figure}[ht]
\centering
\subfigure[]{\includegraphics[width=0.5\textwidth]{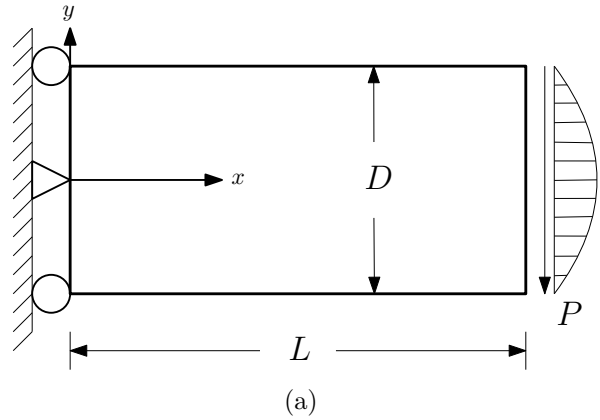}}
\subfigure[]{\includegraphics[width=0.5\textwidth]{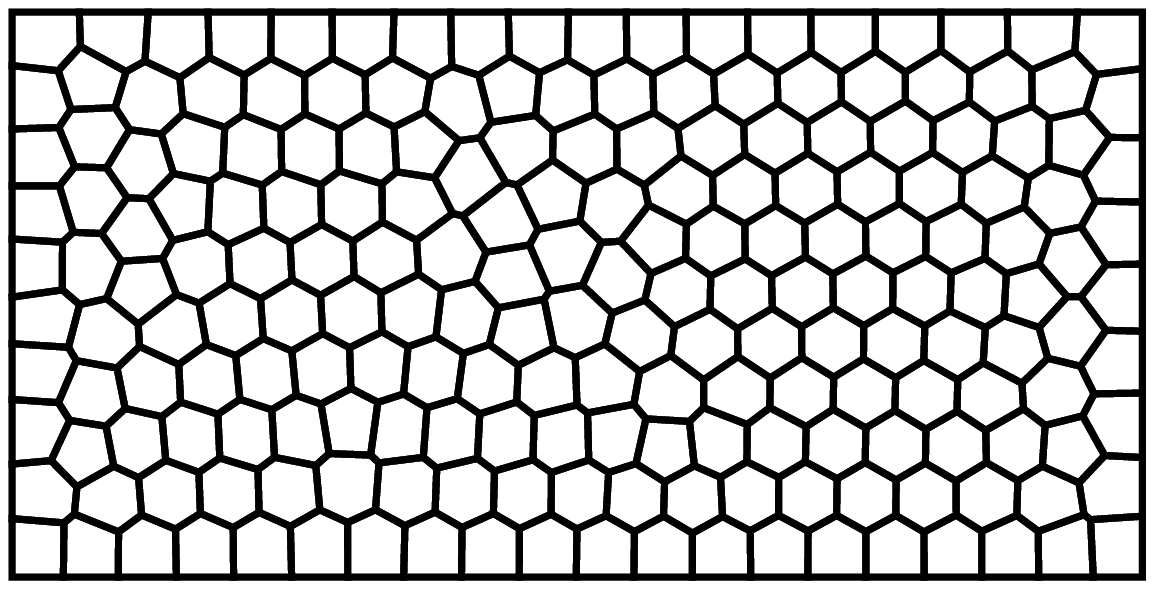}}
\caption{ (a) Cantilever beam: Geometry and boundary conditions and (b) discretization of the domain with polygonal elements.}
\label{fig:cant}
\end{figure}

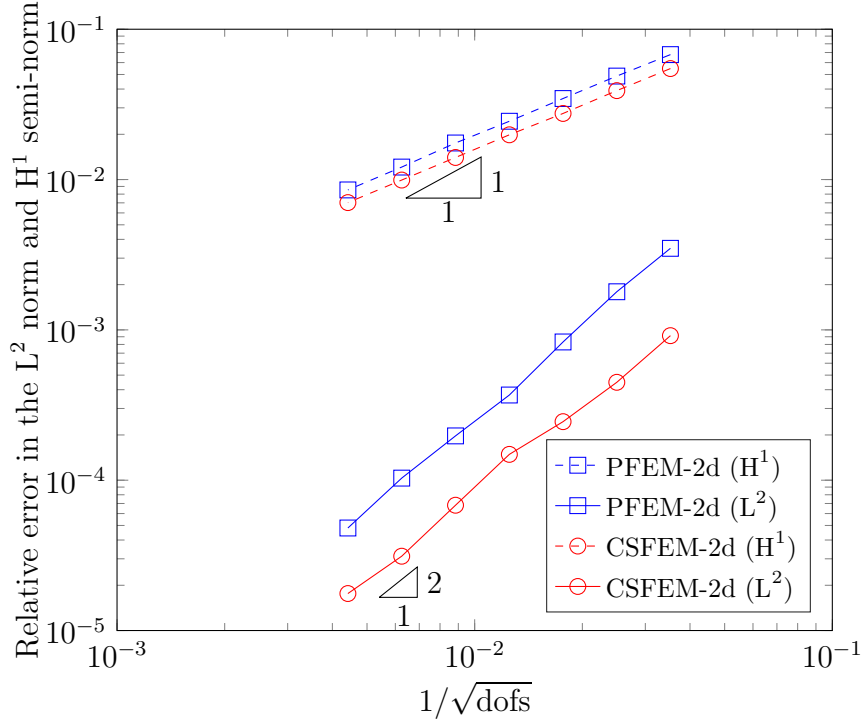
\begin{figure}[h!]
\centering
\newlength\figureheight 
\newlength\figurewidth 
\setlength\figureheight{8cm} 
\setlength\figurewidth{10cm}
\input{./figures/fine_dof_cant_combine_CSFEM_PFEM.tikz}
\caption{Convergence rate for the relative error in the $\textup{L}^2$ norm and the $\textup{H}^1$ semi-norm with mesh refinement for a two-dimensional cantilever beam subjected to end shear}.
\label{fig:cantes}
\end{figure}

Consider a cantilever beam of length $L=$ 8 m, height $D=$ 4 m subjected to a parabolic shear load, $P=$ 250 N at the free end (see \fref{fig:cant}(a)). The beam has the following material properties: Young's modulus, $E=$ 3 $\times$ 10$^7$ Pa and Poisson's ratio, $\nu=$ 0.3. 
The exact solution for this problem is given by \cite{barber2010}:
\begin{align}
u(x,y) &= \frac{P y}{6 EI} \left[ (6L-3x)x + (2+\nu) \left( y^2 - \frac{D^2}{4} \right) \right] \nonumber \\
v(x,y) &= -\frac{P}{6 EI} \left[ 3 \nu y^2(L-x) + (4+5\nu) \frac{D^2x}{4} + (3L-x)x^2 \right] 
\label{eqn:cantisolution}
\end{align}
where $I = D^3/12$ is the moment of inertia and a state of plane stress is considered. The domain is discretized with arbitrary polygonal elements. A representative finite element mesh is shown in \fref{fig:cant}(b). The convergence of the relative error in the displacement and the energy norm with mesh refinement are shown in \fref{fig:cantes} for both the approaches. It is observed that all the methods converge to analytical solution with mesh refinement with optimal convergence rate. 

\subsection{Infinite plate with a circular hole subjected to far field tension}
\begin{figure}[ht]
\centering
\subfigure[]{\includegraphics[width=0.47\textwidth]{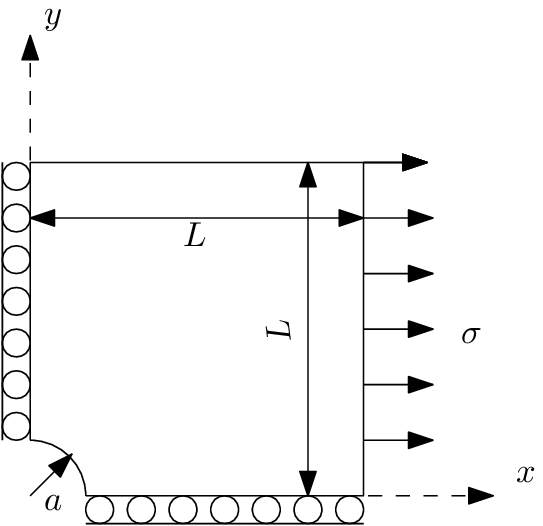}} 
\subfigure[]{\includegraphics[width=0.47\textwidth]{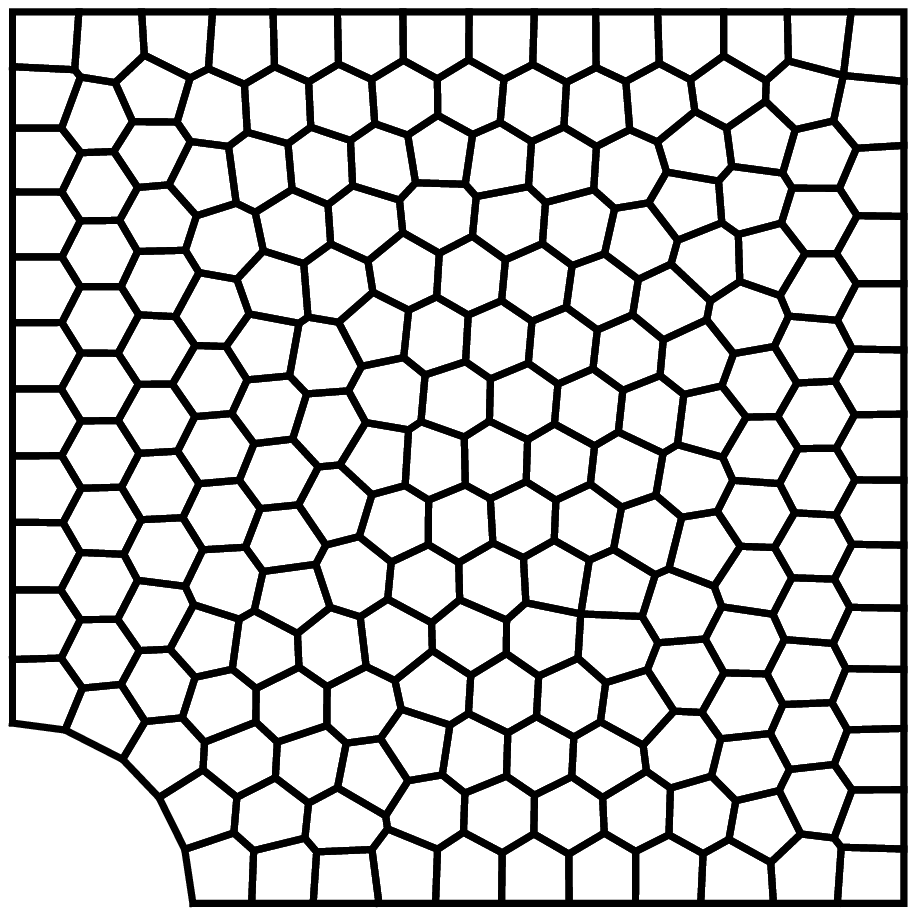}}
\caption{ (a) Infinite plate with a hole: Geometry and boundary conditions and (b) representative finite element discretization: 200 polygonal elements.}
\label{fig:pholegb}
\end{figure}

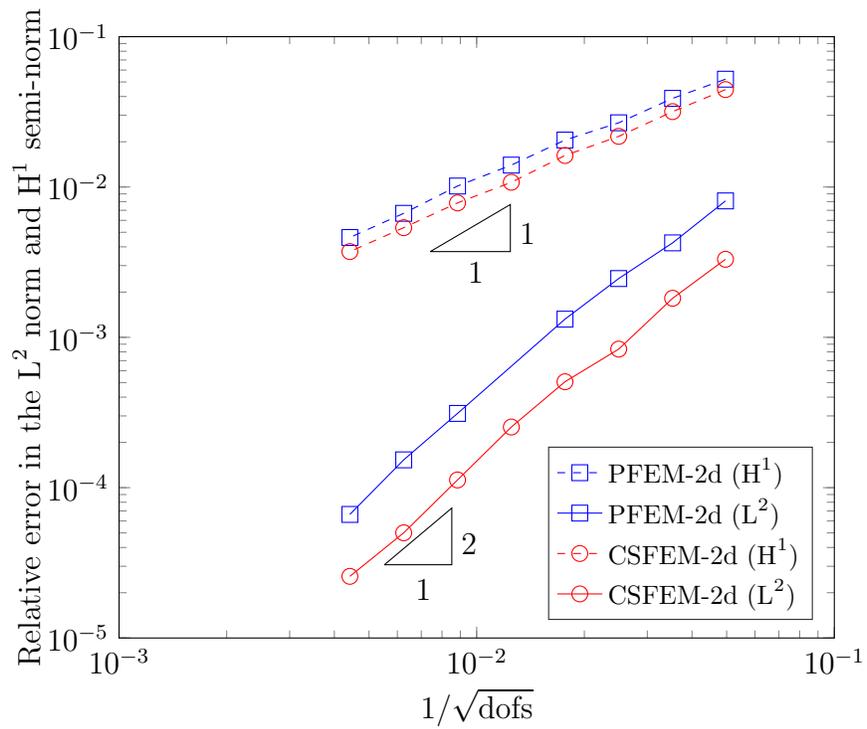
\begin{figure}
\centering
\setlength\figureheight{8cm} 
\setlength\figurewidth{10cm}
\input{./figures/fine_dof_ph_combine_CSFEM_PFEM.tikz}
\caption{Convergence of the relative error in the L$^2$ norm and the H$^1$ seminorm with mesh refinement for a plate with a hole.}
\label{fig:pholeEs}
\end{figure}

Next, consider an infinite plate with a traction free hole subjected to a far field uniformly distributed normal stress $(\sigma=$ 1 Pa$)$ acting along the $x$-axis. A state of plane stress is assumed and the plate has the following material properties: Young's modulus, $E=$ 1000 Pa, Poisson's ratio, $\nu=$ 0.3. The exact solution of the stresses in polar coordinates are given by \cite{barber2010}:
\begin{align}
\sigma_{11}(r,\theta) &= 1 - \frac{a^2}{r^2} \left( \frac{3}{2} (\cos 2\theta + \cos 4\theta) \right) + \frac{3a^4}{2r^4} \cos 4\theta \nonumber \\
\sigma_{22}(r,\theta) &= -\frac{a^2}{r^2} \left( \frac{1}{2}(\cos 2\theta - \cos 4\theta) \right) - \frac{3a^4}{2r^4} \cos 4\theta \nonumber \\
\sigma_{12}(r,\theta) &= -\frac{a^2}{r^2} \left( \frac{1}{2}(\sin 2\theta + \sin 4\theta) \right) + \frac{3a^4}{2r^4} \sin 4\theta
\end{align} 
and the displacements are:
\begin{align}
u_r &= \frac{r}{2E}\left([1+\frac{a^2}{r^2}]+[1-\frac{a^4}{r^4}+4\frac{a^2}{r^2}]\cos2\theta+\nu[1-\frac{a^2}{4^2}]-\nu[1-\frac{a^4}{r^4}]\cos2\theta\right), \nonumber \\
u_\theta &= \frac{r}{2E}\left([1+\frac{a^4}{r^4}+2\frac{a^2}{r^2}]+\nu[1+\frac{a^4}{r^4}-2\frac{a^2}{4^2}]\sin2\theta\right),
\end{align}
where $a$ is the radius of the hole. Due to symmetry, only one quarter of the geometry is modelled as shown in \fref{fig:pholegb}(a). One of the FE meshes is shown in \fref{fig:pholegb}(b), i.e., the discretization with polygonal elements. Analytical solutions are applied on the boundary of the domain and the convergence of the displacements and the stresses are studied.
\fref{fig:pholeEs} shows the convergence of the relative error in both the displacement and the energy norm with reducing mesh size. It is opined that with mesh refinement all the techniques asymptotically converge to reference solution and that the CSFEM yield slightly accurate results. It is noted that the CSFEM requires fewer integration points when compared to the conventional polygonal FEM.

\subsection{Cube subjected to body load}
\begin{figure}[htpb]
\centering
\includegraphics[width=0.4\textwidth]{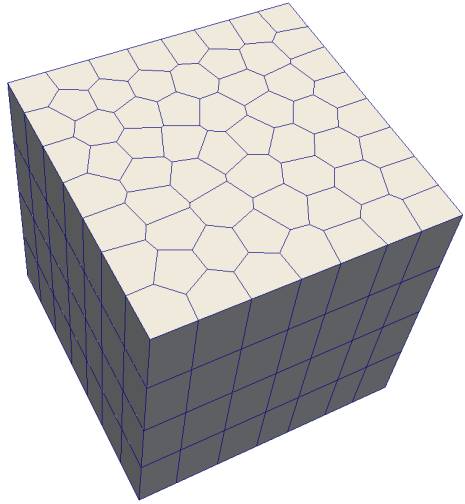}
\caption{Representative polyhedra discretization of a unit cube.}
\label{fig:threeDPT}
\end{figure}

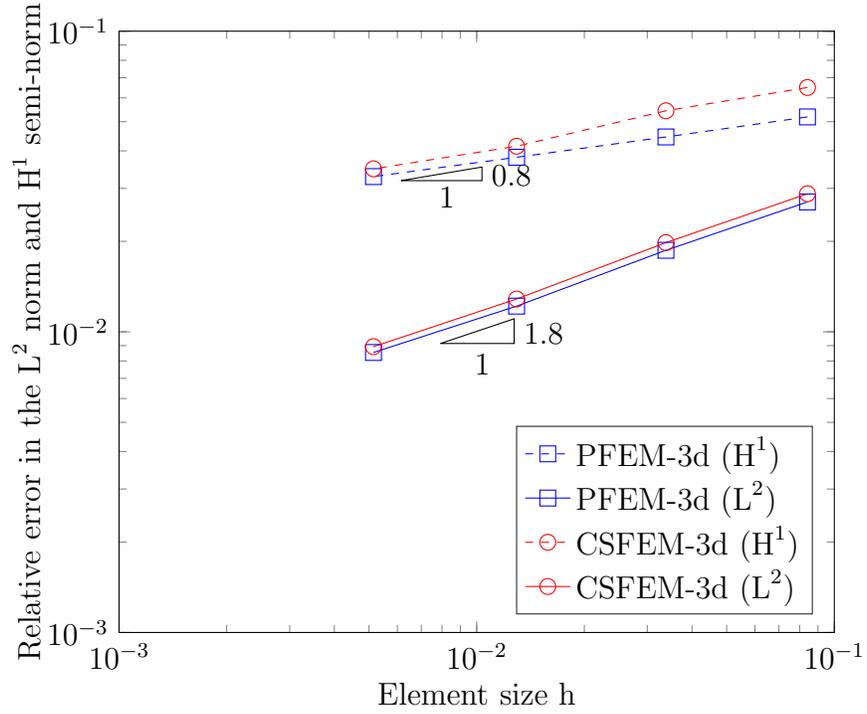
\begin{figure}[h!]
\centering 
\setlength\figureheight{8cm} 
\setlength\figurewidth{10cm}
\input{./figures/QuadPT.tikz}
\caption{Convergence results for the unit cube with polyhedral discretization solved using CSFEM and PFEM. The optimal convergence rates in both the L$^2$ norm and the H$^1$ semi-norm can be seen.}
\label{fig:quadpatch3dres}
\end{figure}

Next, the accuracy and the convergence properties of the polyhedral finite element using CSFEM and the conventional polyhedral finite element method is compared and demonstrated. Consider a cube occupying $\bar{\Omega} =\big\{ (x,y,z) : 0 \leq x \leq 1, -1 \leq y \leq 1 , 0 \leq z \leq 1\big \}$. The boundary of the cube is subjected to the following displacements:
\begin{equation}
\begin{pmatrix} \hat{u} \\ \hat{v} \\ \hat{w} \end{pmatrix} = \begin{pmatrix} 0.1+0.2x+0.2y+0.1z+0.15x^2+0.2y^2+0.1z^2+0.15xy+0.1yz+0.1zx\\ 0.15+0.1x+0.1y+0.2z+0.2x62+0.15y^2+0.1z^2+0.2xy+0.1yz+0.2zx\\ 0.15+0.15x+0.2y+0.1z+0.15x^2+0.1y^2+0.2z^2+0.1xy+0.2yz+0.15zx \end{pmatrix}
\end{equation}
and the domain is subjected to the following body force:
\begin{equation}
\textbf{b} = \begin{pmatrix} -0.3\textbf{C}(1,1)-0.2\textbf{C}(1,2)-0.15\textbf{C}(1,3)-0.6\textbf{C}(4,4)-0.35\textbf{C}(6,6)\\ 
-0.15\textbf{C}(1,2)-0.3\textbf{C}(2,2)-0.2\textbf{C}(2,3)-0.55\textbf{C}(4,4)-0.4\textbf{C}(5,5)\\ 
0.1\textbf{C}(1,3)-0.1\textbf{C}(2,3)-0.4\textbf{C}(3,3)-0.3\textbf{C}(5,5)-0.4\textbf{C}(6,6) \end{pmatrix}
\end{equation}
where \textbf{C} is the material matrix.
The exact solution to \eref{eqn:problem_strong_form} is $\vm{u} =\hat{\vm{u}}$ in the presence of the prescribed body forces. The domain is discretized with polyhedral finite elements for in this study as shown in \fref{fig:threeDPT}. \fref{fig:quadpatch3dres} shows the convergence of the error in the L$^2$ norm and H$^1$ seminorm. As expected, both the polygonal FEM and the CSFEM converges at near optimal rate. It is noted that a very high order numerical integration scheme (i.e. 768 integration points) is employed in case of the conventional polygonal FEM, whilst `only' one integration point per subcell is employed for the CSFEM.

\subsection{Three-dimensional cantilever beam under shear end load}
\begin{figure}[h!]
\centering
\subfigure[]{\includegraphics[scale=0.8]{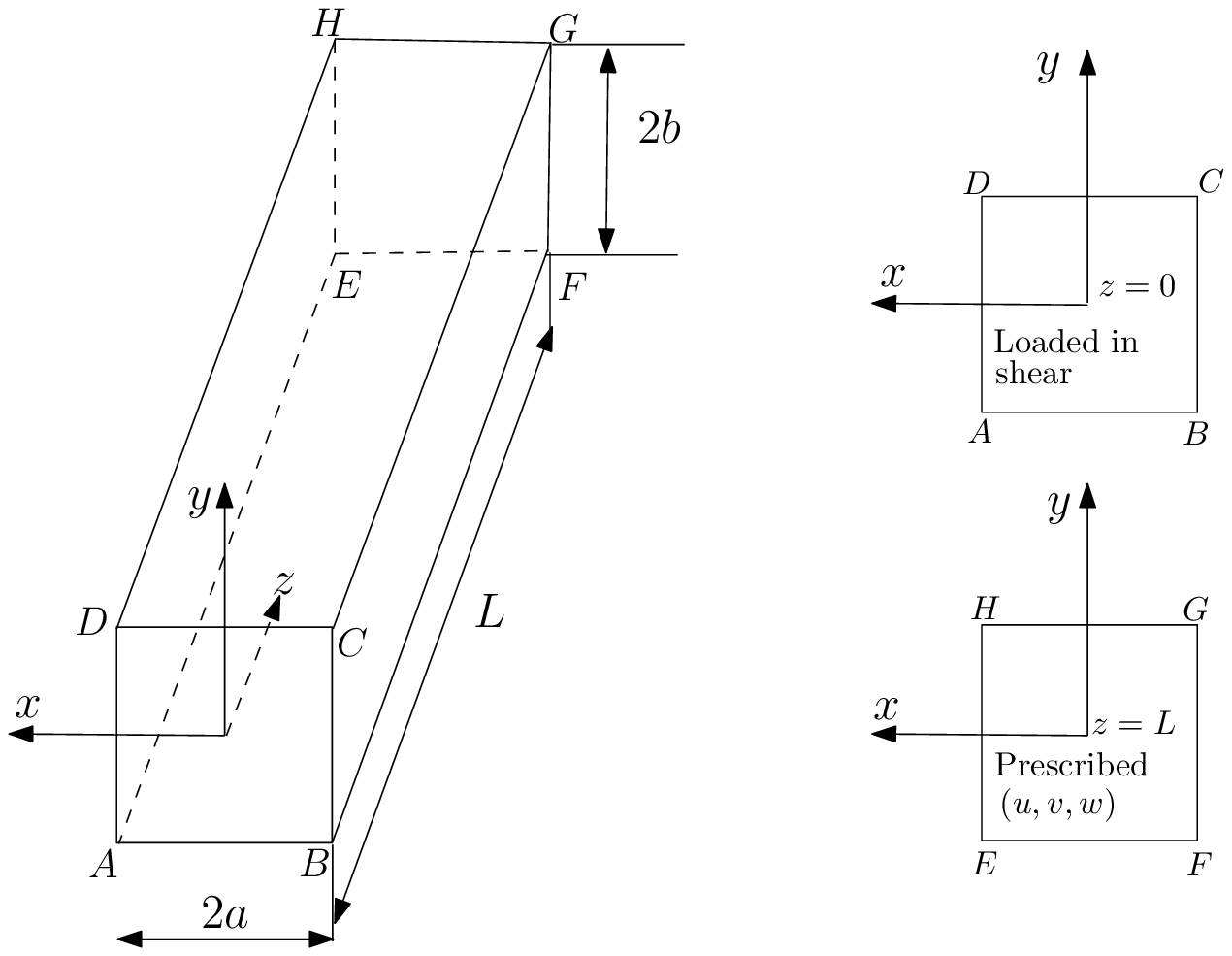}\label{fig:cantprob3d}}
\subfigure[]{\includegraphics[scale=0.35]{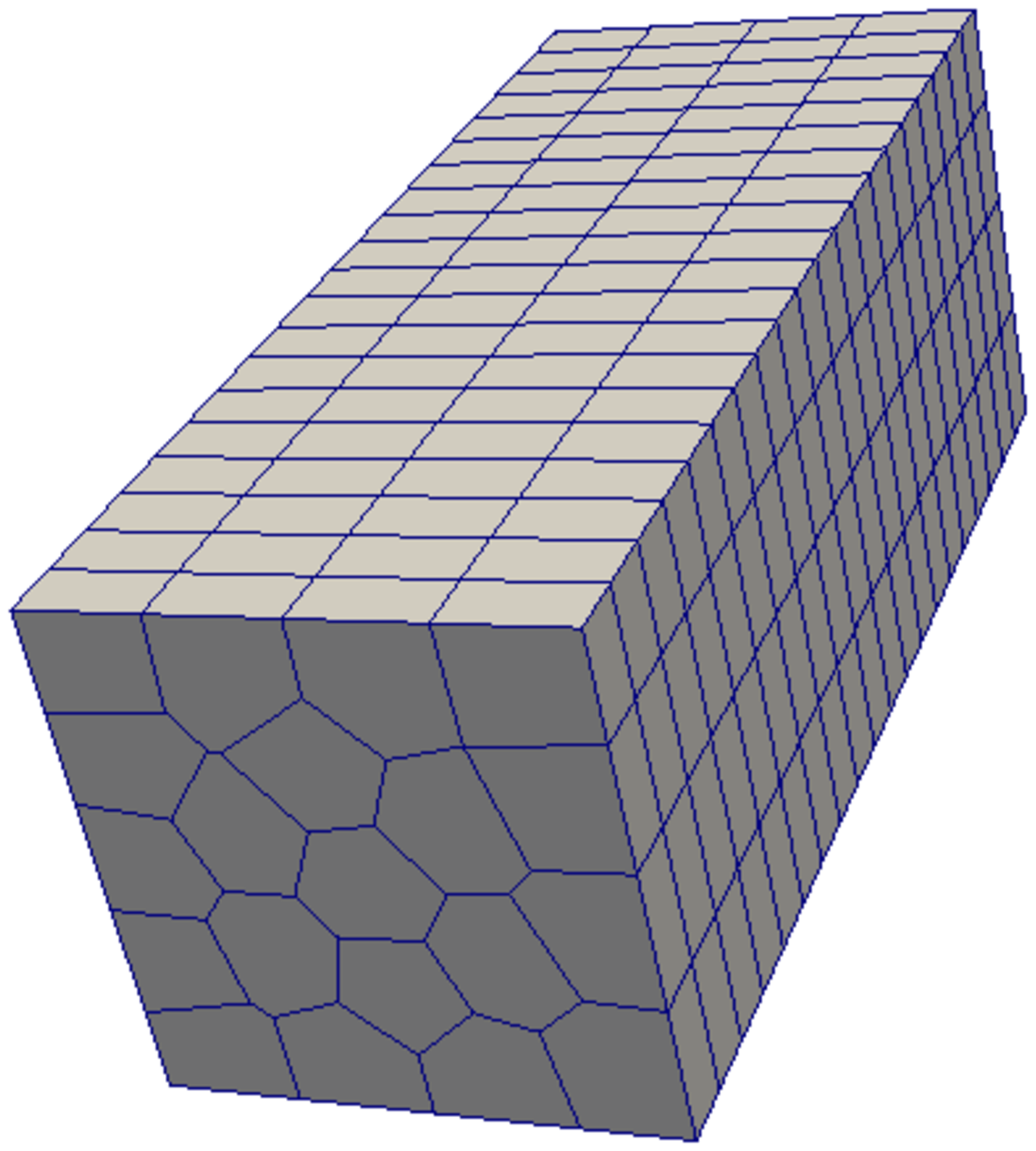}}
\caption{Three-dimensional cantilever beam problem: (a) Geometry and boundary conditions and (b) representative polyhedral mesh used in this study.}
\label{fig:threedcantilever}
\end{figure}

In this example, consider a three-dimensional cantilever beam under shear load at the free end. \fref{fig:cantprob3d} presents the schematic view of the problem and a representative polyhedral discretization of the domain. The domain $\Omega$ for this problem is $[-1,1] \times [-1,1] \times [0,L]$. The material is assumed to be isotropic with Young's modulus, $E=$ 1 N/m$^2$ and Poisson's ratio $\nu=$ 0.3. The beam is subjected to a shear force $F$ at $z=$ 0 and at any cross section of the beam, we have:
\begin{equation}
\int\limits_{-a}^b \int\limits_{-a}^b  \sigma_{yz} ~\rmd x \rmd y = F, \quad \int\limits_{-a}^b \int\limits_{-a}^b \sigma_{zz} y~\rmd x \rmd y = F z.
\end{equation} 
The Cauchy stress field is given by~\cite{barber2010}:
\begin{align}
\sigma_{xx}(x,y,z) &= \sigma_{xy}(x,y,z) = \sigma_{yy}(x,y,z) = 0; \hspace{10pt}
\sigma_{zz}(x,y,z) = \frac{F}{I}yz; \nonumber \hspace{10pt} \\
\sigma_{xz}(x,y,z) &= \frac{2a^2 \nu F}{\pi^2 I (1+\nu)} \sum\limits_{n=0}^\infty \frac{ (-1)^n}{n^2} \sin \left(\frac{n \pi x}{a} \right) \frac{ \sinh \left( \frac{n \pi y}{a} \right)}{ \cosh \left( \frac{n \pi b}{a} \right)} \nonumber \\
\sigma_{yz}(x,y,z) &= \frac{(b^2-y^2)F}{2I} + \frac{ \nu F}{I(1+\nu)} \left[ \frac{ 3x^2-a^2}{6} - \frac{2a^2}{\pi^2} \sum \limits_{n=1}^\infty \frac{ (-1)^n}{n^2} \cos \left(\frac{n \pi x}{a} \right) \frac{ \cosh \left( \frac{n \pi y}{a} \right)}{ \cosh \left( \frac{n \pi b}{a} \right)} \right].
\end{align}

The corresponding displacement field~\cite{bishop2014}:
\begin{align}
u(x,y,z) &= - \frac{\nu F}{EI} xyz; \nonumber \\ 
v(x,y,z) &= \frac{F}{EI} \left[ \frac{\nu (x^2-y^2)z}{2} - \frac{z^3}{6} \right]; \nonumber \\ 
\begin{split}
\hspace{-1.5cm}w(x,y,z) = \frac{F}{EI} \left[ \frac{y (\nu x^2+z^2)}{2} + \frac{\nu y^3}{6} + (1+\nu) \left( b^2y-\frac{y^3}{3} \right) - \frac{\nu a^2 y}{3} - \right. \\ \left. \frac{4 \nu a^3}{\pi^3} \sum\limits_{n=1}^\infty \frac{ (-1)^n}{n^3} \cos \left(\frac{n \pi x}{a} \right) \frac{ \sinh \left( \frac{n \pi y}{a} \right)}{ \cosh \left( \frac{n \pi b}{a} \right)} \right].
\end{split}
\label{eqn:3dcantianadisp}
\end{align}

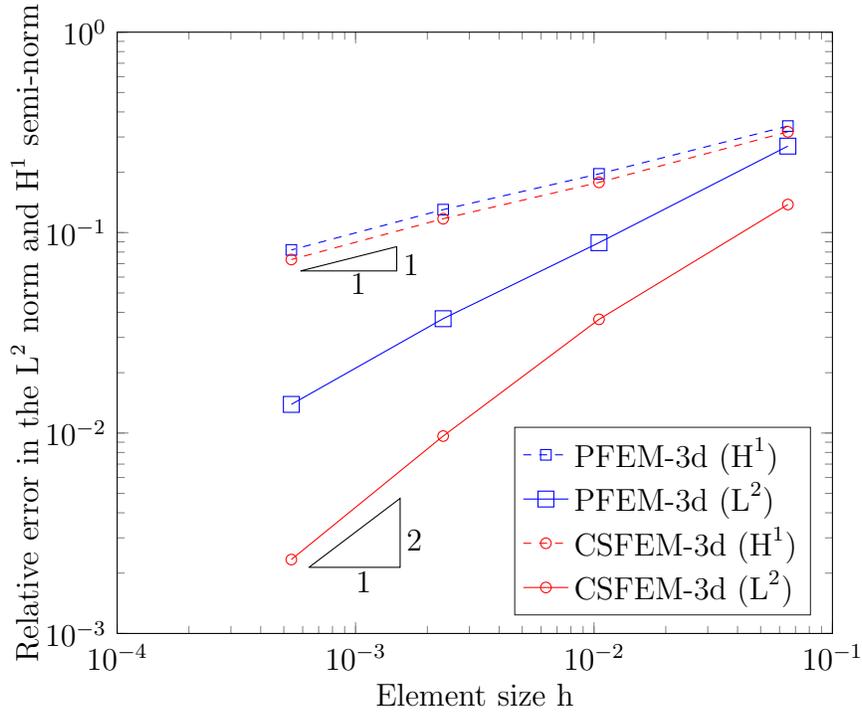
\begin{figure}[htpb]
\centering
\setlength\figureheight{8cm} 
\setlength\figurewidth{10cm}
\input{./figures/Cant3dL2H1.tikz}
\caption{Convergence results for the three-dimensional cantilever beam problem.}
\label{fig:cbeam3dL2}
\end{figure}

where $E$ is the Young's modulus, $\nu$ Poisson's ratio and $I= 4ab^3/3$ is the second moment of area about the $x$-axis. The length of the beam is assumed to be $L=$ 5 m and the shear load is taken as $F=$ 1 N. Analytical displacements given by \eref{eqn:3dcantianadisp} are applied on the beam face at $z=L$ and the beam is loaded in shear on its face at $z=0$. All other faces are assumed to be traction free. \fref{fig:cbeam3dL2} shows the relative error in the L$^2$ norm and H$^1$ seminorm with mesh refinement. It can be seen that both the methods converges optimally. 
\section{Conclusion}
In this work, the CSFEM for two- and three-dimensional linear elastostatic problems is implemented within the commercial finite element software Abaqus. The implementation was based on the user element subroutine feature of Abaqus. The work focused on the main procedures to interact with Abaqus, the structure of the input file and the user element subroutine to compute the smoothed strain-displacement matrix and the corresponding smoothed stiffness matrix. 

The implementation is validated by solving a few benchmark problems in linear elasticity. The UEL implementation has demonstrated the convergence rates and the accuracy levels of the CSFEM as reported in the literature. Although, the paper discusses the results obtained using the CSFEM over polygons in two dimensions and polyhedra (with lateral faces as quadrilateral) in three dimensions, it is noted that the implementation can handle arbitrary polyhedra (in three dimensions), hybrid meshes, i.e., discretizations involving triangles/quadrilaterals/polygonal and tetrahedral/hexhedral/polyhedral elements in two and three dimensions, respectively. 

The developed routines are in the open-source and can be downloaded from\\
\href{Github repository link}{https://github.com/nsundar/SFEM\_in\_Abaqus}. 


\clearpage
\newpage
\bibliographystyle{model1-num-names}
\bibliography{sfem}

\end{document}

%% file: figures/fine_dof_cant_combine_CSFEM_PFEM.tikz
%
%
\begin{tikzpicture}

\begin{axis}[%
width=0.95092\figurewidth,
height=\figureheight,
at={(0\figurewidth,0\figureheight)},
legend style={font=\footnotesize},
scale only axis,
xmode=log,
xmin=0.001,
xmax=0.1,
xminorticks=true,
xlabel={1/$\sqrt{\textup{dofs}}$},
ymode=log,
ymin=1e-05,
ymax=0.1,
yminorticks=true,
ylabel={Relative error in the $\text{L}^\text{2}$ norm and $\mathrm{H}^\text{1}$ semi-norm},
legend style={at={(0.97,0.03)},anchor=south east,legend cell align=left,align=left,draw=white!15!black}
]
\addplot [color=blue,dashed,mark=square,mark options={solid},mark size=3pt]
  table[row sep=crcr]{%
0.0352672807929299	0.067794335536281\\
0.0249688084719461	0.04879083290044\\
0.0176666313334393	0.034599699713728\\
0.0125019535828829	0.024423448536685\\
0.00884021615653596	0.01751011998364\\
0.00625170968564522	0.012107028433053\\
0.0044210583140509	0.008532087455798\\
};
\addlegendentry{$\text{PFEM-2d (H}^\text{1}\text{)}$};

\addplot [color=blue,solid,mark=square,mark options={solid},mark size=3pt]
  table[row sep=crcr]{%
0.0352672807929299	0.003487183290002\\
0.0249688084719461	0.001793734387895\\
0.0176666313334393	0.000831125597662133\\
0.0125019535828829	0.000368717640960067\\
0.00884021615653596	0.000197166001927945\\
0.00625170968564522	0.000103252315859441\\
0.0044210583140509	4.808767040186e-05\\
};
\addlegendentry{$\text{PFEM-2d (L}^\text{2}\text{)}$};

\addplot [color=red,dashed,mark=o,mark options={solid},mark size=3pt]
  table[row sep=crcr]{%
0.0352672807929299	0.054750485511799\\
0.0249688084719461	0.039054787854736\\
0.0176666313334393	0.027523745338987\\
0.0125019535828829	0.019848507025808\\
0.00884021615653596	0.014024037214394\\
0.00625170968564522	0.009924680072858\\
0.0044210583140509	0.007026224379451\\
};
\addlegendentry{$\text{CSFEM-2d (H}^\text{1}\text{)}$};

\addplot [color=red,solid,mark=o,mark options={solid},mark size=3pt]
  table[row sep=crcr]{%
0.0352672807929299	0.000915336623570215\\
0.0249688084719461	0.000447980562120425\\
0.0176666313334393	0.00024514695369188\\
0.0125019535828829	0.000148556833391887\\
0.00884021615653596	6.81648636937754e-05\\
0.00625170968564522	3.12206581566979e-05\\
0.0044210583140509	1.76235675407802e-05\\
};
\addlegendentry{$\text{CSFEM-2d (L}^\text{2}\text{)}$};

\addplot [color=black,solid,forget plot]
  table[row sep=crcr]{%
0.0064210583140509	0.007526224379451\\
0.0104210583140509	0.007526224379451\\
0.0104210583140509	0.0141245277915704\\
0.0064210583140509	0.007526224379451\\
};
\node[right, align=left, inner sep=0mm, text=black]
at (axis cs:0.0080210583140509,0.00609624174735531,0) {$1$};
\node[right, align=left, inner sep=0mm, text=black]
at (axis cs:0.011046321812894,0.00969922865507721,0) {$1$};
\addplot [color=black,solid,forget plot]
  table[row sep=crcr]{%
0.0054110583140509	1.66235675407802e-05\\
0.0069110583140509	1.66235675407802e-05\\
0.0069110583140509	2.64620011717932e-05\\
0.0054110583140509	1.66235675407802e-05\\
};
\node[right, align=left, inner sep=0mm, text=black]
at (axis cs:0.0060110583140509,1.2465089708032e-05,0) {$1$};
\node[right, align=left, inner sep=0mm, text=black]
at (axis cs:0.00732572181289396,2.0694644644133e-05,0) {$2$};
\end{axis}
\end{tikzpicture}%

%% file: figures/fine_dof_ph_combine_CSFEM_PFEM.tikz
%
%
\begin{tikzpicture}

\begin{axis}[%
width=0.95092\figurewidth,
height=\figureheight,
at={(0\figurewidth,0\figureheight)},
scale only axis,
legend style={font=\footnotesize},
xmode=log,
xmin=0.001,
xmax=0.1,
xminorticks=true,
xlabel={1/$\sqrt{\textup{dofs}}$},
ymode=log,
ymin=1e-05,
ymax=0.1,
yminorticks=true,
ylabel={Relative error in the $\textup{L}^\text{2}$ norm and $\textup{H}^\text{1}$ semi-norm},
legend style={at={(0.97,0.03)},anchor=south east,legend cell align=left,align=left,draw=white!15!black}
]
\addplot [color=blue,dashed,mark=square,mark options={solid},mark size=3pt]
  table[row sep=crcr]{%
0.0496291666985465	0.052152804664214\\
0.0353112275773224	0.038890946659227\\
0.0249532564252975	0.026743199406736\\
0.0176611199836457	0.020524339211884\\
0.0125	0.013986139839095\\
0.00884021615653596	0.010149802288197\\
0.00625195404100444	0.006682516283197\\
0.004420453547656	0.004615308524215\\
};
\addlegendentry{$\text{PFEM-2d (H}^\text{1}\text{)}$};

\addplot [color=blue,solid,mark=square,mark options={solid},mark size=3pt]
  table[row sep=crcr]{%
0.0496291666985465	0.008083839827342\\
0.0353112275773224	0.004250282398373\\
0.0249532564252975	0.002460890575788\\
0.0176611199836457	0.001321426116796\\
0.00884021615653596	0.000311843094152946\\
0.00625195404100444	0.000153114402667864\\
0.004420453547656	6.63617314349211e-05\\
};
\addlegendentry{$\text{PFEM-2d (L}^\text{2}\text{)}$};

\addplot [color=red,dashed,mark=o,mark options={solid},mark size=3pt]
  table[row sep=crcr]{%
0.0496291666985465	0.044477900153394\\
0.0353112275773224	0.031736828820404\\
0.0249532564252975	0.021722616053851\\
0.0176611199836457	0.0162060982658\\
0.0125	0.010752843300771\\
0.00884021615653596	0.007838148136284\\
0.00625195404100444	0.005363406105394\\
0.004420453547656	0.00371917898854\\
};
\addlegendentry{$\text{CSFEM-2d (H}^\text{1}\text{)}$};

\addplot [color=red,solid,mark=o,mark options={solid},mark size=3pt]
  table[row sep=crcr]{%
0.0496291666985465	0.003306439800044\\
0.0353112275773224	0.001819331693448\\
0.0249532564252975	0.00083488843098083\\
0.0176611199836457	0.000506493389208646\\
0.0125	0.000253120676047178\\
0.00884021615653596	0.000112519507977243\\
0.00625195404100444	5.00826370417541e-05\\
0.004420453547656	2.57029877225927e-05\\
};
\addlegendentry{$\text{CSFEM-2d (L}^\text{2}\text{)}$};

\addplot [color=black,solid,forget plot]
  table[row sep=crcr]{%
0.007420453547656	0.00371917898854\\
0.012420453547656	0.00371917898854\\
0.012420453547656	0.00764956757551841\\
0.007420453547656	0.00371917898854\\
};
\node[right, align=left, inner sep=0mm, text=black]
at (axis cs:0.009420453547656,0.0027125349807174,0) {$1$};
\node[right, align=left, inner sep=0mm, text=black]
at (axis cs:0.0131656807605154,0.00489603553861975,0) {$1$};
\addplot [color=black,solid,forget plot]
  table[row sep=crcr]{%
0.005520453547656	3.07029877225927e-05\\
0.008520453547656	3.07029877225927e-05\\
0.008520453547656	7.31402676809845e-05\\
0.005520453547656	3.07029877225927e-05\\
};
\node[right, align=left, inner sep=0mm, text=black]
at (axis cs:0.006720453547656,2.088694200553001e-05,0) {$1$};
\node[right, align=left, inner sep=0mm, text=black]
at (axis cs:0.00903168076051536,4.19553365965904e-05,0) {$2$};
\end{axis}
\end{tikzpicture}%

%% file: figures/QuadPT.tikz
%
%
\begin{tikzpicture}

\begin{axis}[%
width=0.95092\figurewidth,
height=\figureheight,
at={(0\figurewidth,0\figureheight)},
scale only axis,
xmode=log,
xmin=0.001,
xmax=0.1,
xminorticks=true,
xlabel={Element size h},
ymode=log,
ymin=0.001,
ymax=0.1,
yminorticks=true,
ylabel={Relative error in the $\text{L}^\text{2}$ norm and $\text{H}^\text{1}$ semi-norm},
legend style={at={(0.97,0.03)},anchor=south east,legend cell align=left,align=left,draw=white!15!black}
]
\addplot [color=blue,dashed,mark=square,mark options={solid},mark size=3pt]
  table[row sep=crcr]{%
0.0841351843463839	0.051860073090518\\
0.0338266544791234	0.044487271755017\\
0.0129331689972914	0.03808506502249\\
0.00514612100579164	0.032839019380365\\
};
\addlegendentry{$\text{PFEM-3d (H}^\text{1}\text{)}$};

\addplot [color=blue,solid,mark=square,mark options={solid},mark size=3pt]
  table[row sep=crcr]{%
0.0841351843463839	0.027043113964162\\
0.0338266544791234	0.018673594120908\\
0.0129331689972914	0.012175707102633\\
0.00514612100579164	0.008538011674638\\
};
\addlegendentry{$\text{PFEM-3d (L}^\text{2}\text{)}$};

\addplot [color=red,dashed,mark=o,mark options={solid},mark size=3pt]
  table[row sep=crcr]{%
0.0841351843463839	0.064972982937546\\
0.0338266544791234	0.054353419209633\\
0.0129331689972914	0.041445531261032\\
0.00514612100579164	0.034849007596466\\
};
\addlegendentry{$\text{CSFEM-3d (H}^\text{1}\text{)}$};

\addplot [color=red,solid,mark=o,mark options={solid},mark size=3pt]
  table[row sep=crcr]{%
0.0841351843463839	0.028795941710463\\
0.0338266544791234	0.019808190061738\\
0.0129331689972914	0.012858301291348\\
0.00514612100579164	0.008927330327625\\
};
\addlegendentry{$\text{CSFEM-3d (L}^\text{2}\text{)}$};

\addplot [color=black,solid,forget plot]
  table[row sep=crcr]{%
0.007927330327625	0.00914612100579164\\
0.012727330327625	0.00914612100579164\\
0.012727330327625	0.0110530233419845\\
0.007927330327625	0.00914612100579164\\
};
\node[right, align=left, inner sep=0mm, text=black]
at (axis cs:0.009847330327625,0.00777420285492289,0) {$1$};
\node[right, align=left, inner sep=0mm, text=black]
at (axis cs:0.0134909701472825,0.00988298922859061,0) {$1.8$};
\addplot [color=black,solid,forget plot]
  table[row sep=crcr]{%
0.00614612100579164	0.031839019380365\\
0.0103461210057916	0.031839019380365\\
0.0103461210057916	0.0353341763887267\\
0.00614612100579164	0.031839019380365\\
};
\node[right, align=left, inner sep=0mm, text=black]
at (axis cs:0.00782612100579164,0.0280183370547212,0) {$1$};
\node[right, align=left, inner sep=0mm, text=black]
at (axis cs:0.0109668882661391,0.033100954251391,0) {$0.8$};
\end{axis}
\end{tikzpicture}%

%% file: figures/Cant3dL2H1.tikz
%
%
\begin{tikzpicture}

\begin{axis}[%
width=0.95092\figurewidth,
height=\figureheight,
at={(0\figurewidth,0\figureheight)},
scale only axis,
xmode=log,
xmin=0.0001,
xmax=0.1,
xminorticks=true,
xlabel={Element size h},
ymode=log,
ymin=0.001,
ymax=1,
yminorticks=true,
ylabel={Relative error in the $\text{L}^\text{2}$ norm and $\text{H}^\text{1}$ semi-norm},
legend style={at={(0.97,0.03)},anchor=south east,legend cell align=left,align=left,draw=white!15!black}
]
\addplot [color=blue,dashed,mark=square,mark options={solid}]
  table[row sep=crcr]{%
0.0650446736340494	0.340087972802895\\
0.0104966116061327	0.196353637900521\\
0.00232917223639432	0.130195388568993\\
0.000537703283956407	0.08195225359175\\
};
\addlegendentry{$\text{PFEM-3d (H}^\text{1}\text{)}$};

\addplot [color=blue,solid,mark=square,mark options={solid},mark size=3pt]
  table[row sep=crcr]{%
0.0650446736340494	0.270063612732915\\
0.0104966116061327	0.089063163701363\\
0.00232917223639432	0.037194739426648\\
0.000537703283956407	0.013896785438453\\
};
\addlegendentry{$\text{PFEM-3d (L}^\text{2}\text{)}$};

\addplot [color=red,dashed,mark=o,mark options={solid}]
  table[row sep=crcr]{%
0.0650446736340494	0.318366199526687\\
0.0104966116061327	0.178097529145969\\
0.00232917223639432	0.117280237207538\\
0.000537703283956407	0.073396358820095\\
};
\addlegendentry{$\text{CSFEM-3d (H}^\text{1}\text{)}$};

\addplot [color=red,solid,mark=o,mark options={solid}]
  table[row sep=crcr]{%
0.0650446736340494	0.138091605732318\\
0.0104966116061327	0.036921837361057\\
0.00232917223639432	0.009656759663967\\
0.000537703283956407	0.002337619519075\\
};
\addlegendentry{$\text{CSFEM-3d (L}^\text{2}\text{)}$};

\addplot [color=black,solid,forget plot]
  table[row sep=crcr]{%
0.000637619519075	0.00213770328395641\\
0.001537619519075	0.00213770328395641\\
0.001537619519075	0.00472069537503238\\
0.000637619519075	0.00213770328395641\\
};
\node[right, align=left, inner sep=0mm, text=black]
at (axis cs:0.000997619519075,0.00179567075852338,0) {$1$};
\node[right, align=left, inner sep=0mm, text=black]
at (axis cs:0.0016298766902195,0.0029244858464618,0) {$2$};
\addplot [color=black,solid,forget plot]
  table[row sep=crcr]{%
0.000587703283956407	0.064396358820095\\
0.00148770328395641	0.064396358820095\\
0.00148770328395641	0.0850880895837738\\
0.000587703283956407	0.064396358820095\\
};
\node[right, align=left, inner sep=0mm, text=black]
at (axis cs:0.000947703283956407,0.0553808685852817,0) {$1$};
\node[right, align=left, inner sep=0mm, text=black]
at (axis cs:0.00157696548099379,0.0686507632992837,0) {$1$};
\end{axis}
\end{tikzpicture}%